\documentclass{article}

\usepackage{graphicx}
\usepackage{amsmath}
\usepackage{amssymb}
\usepackage{color}
\usepackage{hyperref}
\usepackage{epsfig}
\allowdisplaybreaks

\newtheorem{theorem}{Theorem}

\newtheorem{theo}{Theorem}

\newcommand{\NN}{{\mathbb N}}

\newcommand{\ZZ}{{\mathbb Z}}
\newcommand{\RR}{{\mathbb R}}

\newcommand{\bsx}{\boldsymbol{x}}

\newcommand{\notiz}[1]{}

\title{On a multi-dimensional Poissonian pair correlation concept and uniform distribution}
\author{A.\ Hinrichs \footnote{The author is supported by the Austrian Science Fund (FWF) Project F5509-N26, which is a part of the Special Research Program "Quasi-Monte Carlo Methods: Theory and Applications"}, L.\ Kaltenb{\"o}ck \footnote{The author is supported by the Austrian Science Fund (FWF), Project F5507-N26, which is a part of the Special Research Program "Quasi-Monte Carlo Methods: Theory and Applications".}, G.\ Larcher \footnote{The author is supported by the Austrian Science Fund (FWF), Project F5507-N26, which is a part of the Special Research Program "Quasi-Monte Carlo Methods: Theory and Applications" and Project I1751-N26.}, \\  W.\ Stockinger \footnote{The author is supported by the Austrian Science Fund (FWF), Project F5507-N26, which is a part of the Special Research Program "Quasi-Monte Carlo Methods: Theory and Applications".}, M.\ Ullrich}
\date{}
\begin{document}
\maketitle 
\begin{abstract}
The aim of the present article is to introduce a concept which allows to generalise the notion of Poissonian pair correlation, a second-order equidistribution property, to higher dimensions. Roughly speaking, in the one-dimensional setting, the pair correlation statistics measures the distribution of spacings between sequence elements in the unit interval at distances of order of the mean spacing $1/N$. In the $d$-dimensional case, of course, the order of the mean spacing is $1/N^{\frac{1}{d}}$, and --in our concept-- the distance of sequence elements will be measured by the supremum-norm. 

Additionally, we show that, in some sense, almost all sequences satisfy this new concept and we examine the link to uniform distribution. The metrical pair correlation theory is investigated and it is proven that a class of typical low-discrepancy sequences in the high-dimensional unit cube do not have Poissonian pair correlations, which fits the existing results in the one-dimensional case. 
\end{abstract}
\section{Introduction and statement of results}
The concept of Poissonian pair correlations has its origin in quantum mechanics, where the spacings of energy levels of integrable systems were studied. See for example \cite{not1} and the references cited therein for detailed information on that topic. Rudnick and Sarnak first studied this concept from a purely mathematical point of view and over the years the topic has attracted wide attention, see e.g., \cite{not7, not10, not11, not12, not13}. \\ \\

Let $\| \cdot \|$ denote the distance to the nearest integer. A sequence $(x_n)_{n \in \NN}$ of real numbers in $[0,1)$ has Poissonian pair correlations if the pair correlation statistics
\begin{equation}\label{eq:pc1}
F_N(s):= \frac{1}{N} \# \left\lbrace 1  \leq l \neq m \leq N: \| x_l - x_m \| \leq \frac{s}{N} \right\rbrace 
\end{equation} 
tends to $2s$, for every $s \geq 0$, as $N \to \infty$. \\ \\
Let now $d \geq 2$ be an integer denoting the dimension of the problem setting. In the sequel, we indicate by bold symbols that we work with $d$-dimensional vectors of real numbers or random variables. We extend the above notion to sequences $(\boldsymbol{x}_n)_{n \in \NN}$ in the $d$-dimensional unit cube $[0,1)^d$. Subsequently, we denote by $\| \cdot \|_{\infty}$ a supremum-norm of a $d$-dimensional vector, i.e., in our case for some $\boldsymbol{x}= (x_1, \ldots, x_d) \in \mathbb{R}^{d}$,
\begin{equation*}
\| \boldsymbol{x} \|_{\infty} := \max(\| x_1 \|, \ldots, \| x_d \|),
\end{equation*}
where we recall that $\| \cdot \|$ denotes the distance to the nearest integer. 

We say that a sequence $(\boldsymbol{x}_n)_{n \in \NN} \in [0,1)^d$ has Poissonian pair correlations if the multi-dimensional pair correlation statistics
\begin{equation*}
F_N^{(d)}(s):= \frac{1}{N} \# \left\lbrace 1  \leq l \neq m \leq N: \| \boldsymbol{x}_l - \boldsymbol{x}_m \|_{\infty} \leq \frac{s}{N^{1/d}} \right\rbrace 
\end{equation*}
tends to $(2s)^d$, for every $s \geq 0$, as $N \to \infty$. \\ \\ 
In the one-dimensional case, it is well-known that for a sequence of i.i.d. random variables $(X_i)_{i \in \NN}$ having uniform distribution on $[0,1)$, $F_N(s)$ tends to $2s$, as $N \to \infty$, almost surely (see e.g., \cite{not2}). We prove a multi-dimensional analogue.
\begin{theorem}\label{TH:TH1}
Let $(\boldsymbol{X}_i)_{i \in \NN}$ a sequence of i.i.d. random variables having uniform distribution on $[0,1)^d$, then for all $s>0$, we have $F_N^{(d)}(s) \to (2s)^d$, as $N \to \infty$, almost surely.  
\end{theorem}
The notion Poissonian pair correlation has attracted renewed interest in the last few years, due to its connection to several mathematical fields, such as Diophantine approxmation, additive combinatorics and uniform distribution (see e.g., \cite{not1, not3, not5, not7, not8, not4}). The link between the concept of uniform distribution modulo 1 and Poissonian pair correlation has been studied in the one-dimensional case. Due to a result by Grepstad and Larcher \cite{not6} (see also \cite{not2, not14}), we know that a sequence which satisfies that (\ref{eq:pc1}) tends to $2s$, for every $s > 0$, as $N \to \infty$, is also uniformly distributed in $[0,1)$, i.e., it satisfies 
\begin{equation*}
\lim_{N \to \infty} \frac{1}{N} \# \lbrace 1 \leq n \leq N: x_n \in [a,b) \rbrace = b-a 
\end{equation*}
for all $0 \leq a < b \leq 1$. The above presented multi-dimensional concept of Poissonian pair correlation also implies uniform distribution of a sequence in $[0,1)^d$. 
\begin{theorem}\label{TH:TH2}
Let $(\boldsymbol{x}_n)_{n \in \NN} \in [0,1)^d$ be such that for every $s \in \NN$ we have that
\begin{equation*}
\lim_{N \to \infty} F_N^{(d)}(s) = (2s)^d,
\end{equation*}
then $(\boldsymbol{x}_n)_{n \in \NN}$ is uniformly distributed in $[0,1)^d$
\end{theorem}
It turns out by the proof of this theorem that it is sufficient to have the Poissonian property for positive integer-valued $s$ only in order to deduce uniform distribution for a sequence $(\boldsymbol{x}_n)_{n \in \NN} \in [0,1)^d$ (the same holds in the one-dimensional case as well). \\ \\
Classical low-discrepancy sequences in $[0,1)$, e.g., the van der Corput sequence, the Kronecker sequence $(\lbrace n \alpha \rbrace)_{n \in \NN}$ and digital $(t,1)$-sequences in base $b \geq 2$, do not have Poissonian pair correlation (see e.g., \cite{not9}). We will derive an analogous result for the multi-dimensional version of the Kronecker sequence.   
\begin{theorem}\label{TH:TH3}
For any choice of $\alpha_1, \ldots, \alpha_d \in \mathbb{R}$ the $d$-dimensional Kronecker sequence 
\begin{equation*}
(\boldsymbol{x}_n)_{n \in \NN}:=((\lbrace n\alpha_1 \rbrace, \ldots, \lbrace n \alpha_d \rbrace))_{n \in \NN},
\end{equation*} 
where $\lbrace \cdot \rbrace$ denotes the fractional part of a real number, does not have Poissonian pair correlations. 
\end{theorem}
We also strongly believe that other high-dimensional low-discrepancy sequences such as $(t,s)$-sequences and the Halton sequences do not have Poissonian pair correlations. \\ \\
This new concept of course raises several further questions. E.g., in the one-dimensional case, it is known that for almost all choices of $\alpha$ the sequence $( \lbrace f(n) \alpha \rbrace)_{n \in \NN}$, where $f(x)$ is a polynomial of degree at least $2$ with integer coefficients, has Poissonian pair correlations (\cite{not11}). If a $d$-dimensional polynomial $\boldsymbol{p}(x)= ( p_1(x), \ldots, p_d(x))$ ($p_i(x)$ are all real polynomials) has the property that for each lattice point $\boldsymbol{h} \in \ZZ^d$,  $\boldsymbol{h} \neq \boldsymbol{0}$ the polynomial $\left\langle \boldsymbol{h}, \boldsymbol{p}(x) \right\rangle $ has at least one non-constant term with irrational coefficient, then 
\begin{equation*}
((\lbrace p_1(n)\rbrace, \ldots, \lbrace p_d(n) \rbrace))_{n \in \NN}
\end{equation*}
is uniformly distributed in $[0,1)^d$ (see e.g., \cite{not18}). Therefore, in analogy to the one-dimensional case, it would be natural to expect that for an integer polynomial $f(x)$ with degree at least $2$, the sequence 
\begin{equation*}
(( \lbrace f(n)\alpha_1 \rbrace, \ldots, \lbrace f(n) \alpha_d \rbrace ))_{n \in \NN},
\end{equation*}
has Poissonian pair correlation for almost all choices of $\alpha_1, \ldots, \alpha_d$ (this claim will be a consequence of Theorem \ref{TH:THMet}). 

To be more general, let $(a_n)_{n \in \NN}$ be an increasing sequence of distinct integers and $\boldsymbol{\alpha}=(\alpha_1, \ldots, \alpha_d)$. Then we consider a sequence of the form 
\begin{equation*}
( \lbrace a_n \boldsymbol{\alpha} \rbrace)_{n \in \NN}:=((\lbrace a_n\alpha_1 \rbrace, \ldots, \lbrace a_n \alpha_d \rbrace))_{n \in \NN},
\end{equation*}
which is uniformly distributed for almost all choices of $\alpha_1, \ldots, \alpha_d$.

In the one-dimensional case the metrical pair correlation theory of such sequences is strongly linked to the additive energy of a finite integer set $A$, denoted by $E(A)$. The additive energy $E(A)$ is defined as  
\begin{equation*}
E(A):= \sum_{a+b=c+d} 1,
\end{equation*}
where the sum is extended over all quadruples $(a,b,c,d) \in A^4$. This connection was discovered by Aistleitner, Larcher and Lewko, who, roughly speaking, proved in \cite{not5} that if the first $N$ elements of an increasing sequence of distinct integers $(a_n)_{n \in \NN}$, have an arbitrarily small energy saving, then $(\lbrace a_n \alpha \rbrace)_{n \in \NN}$ has Poissonian pair correlations for almost all $\alpha$. Recently, Bloom and Walker (see \cite{not20}) improved over this result by showing the following theorem. 
\begin{theo}\label{TH:THWalk}
There exists an absolute positive constant $C$ such that the following is true. Let $A_N$ denote the first $N$ elements of $(a_n)_{n \in \NN}$ and suppose that 
\begin{equation*}
E(A_N) = \mathcal{O}\left(\frac{N^3}{(\log N)^C} \right),
\end{equation*} 
then for almost all $\alpha$, $(\lbrace a_n \alpha \rbrace)_{n \in \NN}$ has Poissonian pair correlations.
\end{theo}
The proof of this results relies on a new bound for GCD sums with $\alpha=1/2$, 
which improves over the bound by Bondarenko and Seip (see, \cite{not21}), 
if the additive energy of $A_N$ is sufficiently large.
Note that the constant $C$ was not specified in the above mentioned paper, 
but the authors thereof conjecture that Theorem~\ref{TH:THWalk} holds for $C>1$ already. 
This result would be best possible. To see this, consider the sequence $(p_n)_{n\in\NN}$ 
of all primes. It is known that $(p_n)_{n=1}^N$ has additive energy of exact order 
$N^3/(\log N)$, but $(\lbrace p_n \alpha \rbrace)_{n \in \NN}$ is not metric Poissonian, i.e., there exists a set $\Omega$ of full Lebesgue measure, such that for all $\alpha \in \Omega$, $(\lbrace p_n \alpha \rbrace)_{n \in \NN}$ does not have Poissonian pair correlations (see, \cite{not4}).
\medskip

Naturally, we would also expect that under this condition on the additive energy, the sequence 
\begin{equation*}
( \lbrace a_n \boldsymbol{\alpha} \rbrace)_{n \in \NN}
\end{equation*}
has Poissonian pair correlations for almost all instances and, in fact, we have the following result:
\begin{theorem}\label{TH:THMet}
There exists an absolute positive constant $C$ such that the following is true. Let $A_N$ denote the first $N$ elements of $(a_n)_{n \in \NN}$ and suppose that 
\begin{equation*}
E(A_N) = \mathcal{O}\left(\frac{N^3}{(\log N)^C} \right),
\end{equation*} 
then for almost all choices of $\boldsymbol{\alpha}=(\alpha_1, \ldots, \alpha_d) \in \mathbb{R}^d$,
\begin{equation*}
( \lbrace a_n \boldsymbol{\alpha} \rbrace)_{n \in \NN}
\end{equation*}
has Poissonian pair correlations.
\end{theorem}
However, if the additive energy is of maximal order, i.e., if we have $E(A_N)=\Omega(N^3)$, then there is no $\alpha$ such that $(\lbrace a_n \alpha \rbrace)_{n \in \NN}$ has Poissonian pair correlations, see \cite{not8}. The approach used in \cite{not8} can be generalised to arbitrary dimensions.
\begin{theorem}\label{TH:THMAX}
If  $E(A_N)=\Omega(N^3)$, then for any choice of $\boldsymbol{\alpha} =(\alpha_1, \ldots, \alpha_d) \in \mathbb{R}^d$ the sequence
\begin{equation*}
( \lbrace a_n \boldsymbol{\alpha} \rbrace)_{n \in \NN},
\end{equation*} 
does not have Poissonian pair correlations. 
\end{theorem}

Of course, Theorem \ref{TH:TH3} could also immediately be deduced by 
Theorem \ref{TH:THMAX}. However, we also include the explicit proof of 
Theorem \ref{TH:TH3}, since it gives an intuitive feeling for the multi-dimensional 
Poissonian pair correlation concept.

\section{Proof of Theorem \ref{TH:TH1}}
We adopt some of the steps of \cite{not15} (Chapter 2, and in particular Theorem 2.3.). \\ \\
As $(\boldsymbol{X}_i)_{i \in \NN}$ is a sequence of i.i.d.\ random variables having uniform distribution on $[0,1)^d$, we have
\begin{align*}
\mathbb{E}\left(F_N^{(d)}(s)\right) &= \mathbb{E}\left( \frac{1}{N} \# \left\lbrace 1  \leq l \neq m \leq N: \| \boldsymbol{X}_l - \boldsymbol{X}_m \|_{\infty} \leq \frac{s}{N^{1/d}} \right\rbrace \right) \\
&= \frac{1}{N}N(N-1) \int_{[0,1)^d} \int_{x_1 -s/N^{1/d}}^{x_1 + s/N^{1/d}} \ldots \int_{x_d -s/N^{1/d}}^{x_d + s/N^{1/d}} 1 \ d \boldsymbol{y} \ dx_d \ldots \ d x_1 \\ 
&=\frac{N-1}{N}(2s)^d.
\end{align*} 
To compute higher-order moments of the pair correlation statistics, 
we consider a different representation of $F_N^{(d)}(s)$ using a Fourier analytic approach.
For this let $\boldsymbol{I}\colon\RR^d\to\{0,1\}$ be the indicator function of the set $[-1/2,1/2]^d$ and 
note that its Fourier transform is given by 
\[
\mathcal{F}\boldsymbol{I}(\boldsymbol{\xi}) \,=\, \prod_{i=1}^d \frac{\sin(\pi \xi_i)}{\pi \xi_i}, 
\]
where $\boldsymbol{\xi}=(\xi_1,\dots,\xi_d)\in\RR^d$ and, for $f\in L_1(\RR^d)$, 
\[
\mathcal{F}f(\boldsymbol{\xi}) \,:=\, \int_{\RR^d} f(\bsx)\, e\bigl(-\left\langle \boldsymbol{\xi},\bsx \right\rangle\bigr)\, d\bsx,
\]
with $e(x):= \exp(2 \pi i x)$.

We write, using the Poisson summation formula, 
\begin{align*}
F_N^{(d)}(s) \,&=\, \frac{1}{N} \sum_{1 \leq k \neq l \leq N} \sum_{\boldsymbol{q} \in \ZZ^{d}} 
	\boldsymbol{I}\left(\frac{(\boldsymbol{X}_k - \boldsymbol{X}_l + \boldsymbol{q})N^{1/d}}{2s} \right) \\ 
&=\, \frac{(2s)^d}{N^2} \sum_{1 \leq k \neq l \leq N} \sum_{\boldsymbol{r} \in \ZZ^d} \mathcal{F}\boldsymbol{I}\left( \frac{(2s) \boldsymbol{r}}{N^{1/d}} \right) e(\left\langle \boldsymbol{r},\boldsymbol{X}_k -\boldsymbol{X}_l \right\rangle).
\end{align*}
Hence, we get, with 
\begin{align*}
&\boldsymbol{X}_{\circ}=(X^{(1)}_{\circ}, \ldots, X^{(d)}_{\circ}) \quad \text{ and} \\
&\boldsymbol{r}=(r_1,\ldots, r_d)\quad  (\text{analogously for } \boldsymbol{r}{'} ),
\end{align*}
that
\begin{align*}
& \mathbb{E}\left[ \left( F_N^{(d)}(s) -  \mathbb{E}\left( F_N^{(d)}(s) \right)\right)^2 \right] \\
&= \frac{(2s)^{2d}}{N^4} \sum_{\substack{1 \leq k, l, m, n  \leq  N \\ k \neq l, \ m \neq n }} 
			\sum_{\boldsymbol{r}, \boldsymbol{r}{'}  \in \ZZ^{d} \setminus \lbrace \boldsymbol{0} \rbrace}  
			\mathcal{F}\boldsymbol{I}\left(\frac{(2s) \boldsymbol{r}}{N^{1/d}} \right) 
			\mathcal{F}\boldsymbol{I}\left(\frac{(2s) \boldsymbol{r}{'}}{N^{1/d}} \right)  \times \\
&\qquad\quad \times \mathbb{E}\left[e \left(\sum_{i=1}^{d} r_i\left(X^{(i)}_k -X^{(i)}_l \right) - r_i{'}\left(X^{(i)}_m -X^{(i)}_n \right) \right) \right]. 
\end{align*}
Let $ \vartheta \subset \mathcal{D}:=\lbrace 1, \ldots, d\rbrace$ denote the subset of 
indices $i$ for which $r_i = r_i{'}=0$. 
Note that $ 0\leq |\vartheta |<d$.  Moreover, 
\begin{align*}
&\mathbb{E}\left[e \left(\sum_{i=1}^{d} r_i\left(X^{(i)}_k -X^{(i)}_l \right) - r_i{'}\left(X^{(i)}_m -X^{(i)}_n \right) \right) \right] = \\
&=
\begin{cases}
1, \ &\text{if } 0 \neq r_i=r_i{'} , \ k=m, \ l=n, \\ 
 \ &\text{or  if } 0 \neq r_i=-r_i{'}, \ k=n, \ l=m, \\
 \ &\text{for } i \notin \vartheta \\
0, \ &\text{otherwise,}
\end{cases}
\end{align*}
and, clearly, 
\[
\bigl|\mathcal{F}\boldsymbol{I}(\boldsymbol{\xi})\bigr| 
	\,\le\, \prod_{i=1}^d \min\left(1,\frac{1}{\pi |\xi_i|}\right).  
\]
Using rather simple computations, we obtain 
\begin{align*}
\mathbb{E}\left[ \left(F_N^{(d)}(s) -  \mathbb{E}\left( F_N^{(d)}(s) \right)\right)^2 \right] 
\,&=\, \frac{(2s)^{2d}}{N^4}\sum_{\vartheta \subset \mathcal{D}}  \mathcal{O}\left( N^2 \left(\frac{N^{1/d}}{2s}\right)^{d-|\vartheta|} \right) \\
&=\, \mathcal{O}\left(\frac{\max(s^d, s^{2d-1})}{N} \right), 
\end{align*}
where the implied constant depends on the dimension $d$, but is independent of $s$ and $N$. 

Using Chebyshev's inequality, we obtain that there exists a constant $c>0$ 
such that for all $\epsilon>0$, $N$, $s$, 
\begin{equation*}
\mathbb{P}\left(\left|F_N^{(d)}(s) - (2s)^d \right| \geq \epsilon  \right) 
\,\leq\, c\, \frac{\max(s^d, s^{2d-1})}{\epsilon^2 N}.
\end{equation*}
To prove now almost sure convergence, one can apply the arguments used in \cite{not5, not13}. We fix a $\gamma >0$ and define a subsequence $N_M$ along the integers, for $M \geq 1$, as
\begin{equation*}
N_{M}:=M^{1+\gamma}.
\end{equation*} 
The variance estimate from above combined with Chebyshev's inequality and the first Borel-Cantelli lemma allow to deduce that, for all $s >0$,
\begin{equation*}
\lim_{M \to \infty} F_{N_M}^{(d)}(s) = (2s)^d, \text{ almost surely}, 
\end{equation*}
i.e., we have now almost sure convergence along a subsequence of the integers. 
For $N$, with $N_M \leq N \leq  N_{M+1}$, we use the trivial bounds 
\begin{align*}
N_M F_{N_M}^{(d)}\left(\frac{N_M}{N_{M+1}}s\right) &\leq N F_{N}^{(d)}\left(s\right) \\
& \leq N_{M+1} F_{N_{M+1}}^{(d)}\left(\frac{N_{M+1}}{N_M}s\right).
\end{align*}
Since $N_{M+1}/N_M \to 1$, as $M \to \infty$, we also get 
\begin{equation*}
\lim_{N \to \infty} F_{N}^{(d)}(s) = (2s)^d, \text{ almost surely}.
\end{equation*}
\hfill $\square$

\section{Proof of Theorem \ref{TH:TH2}}
We carry out the proof for $d=2$, for an arbitrary $d$ the arguments run quite analogously.  
To prove the theorem, we assume in the contrary that $(\boldsymbol{x}_n)_{n \in \NN}$ is not uniformly distributed and will derive a contradiction. Due to this assumption, there exists an $\epsilon > 0$ and $\alpha, \beta$ with $0 < \alpha, \beta <1$ such that 
\begin{equation*}
\left| \frac{1}{N}  \# \lbrace 1 \leq n \leq N : \ \boldsymbol{x}_n \in [0, \alpha) \times [0, \beta) - \alpha \beta \rbrace \right| > \epsilon,
\end{equation*}   
for infinitely many $N$. Hence, we can assume that for an increasing sequence of integers $(N_i)_{i \in \NN}$ we have 
\begin{equation*}
\frac{1}{N_i}  \# \lbrace 1 \leq n \leq N_i : \ \boldsymbol{x}_n \in [0, \alpha) \times [0, \beta) \rbrace  \leq \alpha \beta - \epsilon
\end{equation*}
(The case that we have "$\geq \alpha \beta + \epsilon$" in the above expression can be treated analogously.) Let $N:=N_i$ for some $i \geq 1$ and assume for simplicity that $\sqrt{N}$ is an integer. For $0 \leq i,j < \sqrt{N}$ let 
\begin{equation*}
A_{i,j}:= \# \left\lbrace 1 \leq n \leq N : \ \boldsymbol{x}_n \in \left[ \frac{i}{\sqrt{N}}, \frac{i+1}{\sqrt{N}} \right) \times \left[ \frac{j}{\sqrt{N}}, \frac{j+1}{\sqrt{N}} \right) \right\rbrace.
\end{equation*}
If $i$ and/or $j \geq \sqrt{N}$, we set 
\begin{equation*}
A_{i,j} := A_{i \text{ mod } \sqrt{N}, \ j \text{ mod } \sqrt{N}}.
\end{equation*}  
Then for all integers $s$ and $N$ large enough, we have
\begin{equation*}
N F_N^{(d)}(s) \geq \sum_{i,j =0}^{\sqrt{N}-1} A_{i,j} (A_{i,j}-1) + 2A_{i,j}B_{i,j}, \text{ where}
\end{equation*} 
\begin{equation*}
B_{i,j}:= \sum_{u=1}^{s-1} \sum_{v=-(s-1)}^{s-1} A_{i+u,j+v} + \sum_{v=1}^{s-1}A_{i,j+v}.
\end{equation*}
Hence, 
\begin{equation*}
N F_N^{(d)}(s) \geq \sum_{i,j=0}^{\sqrt{N}-1}\left( (A_{i,j}+B_{i,j})^2 - B_{i,j}^2 \right) -N,
\end{equation*}
as we have 
\begin{equation}\label{eq:eq1}
\sum_{i,j=0}^{\sqrt{N}-1} A_{i,j} = N.
\end{equation}
If $\alpha = \frac{a}{\sqrt{N}}$ and $\beta= \frac{b}{\sqrt{N}}$ (assume for simplicity that $a$ and $b$ are positive integers), then we have 
\begin{equation}\label{eq:eq2}
\sum_{i=0}^{a-1} \sum_{j=0}^{b-1} A_{i,j} \leq N(\alpha \beta - \epsilon)
\end{equation}
and 
\begin{equation}\label{eq:eq3}
\sum_{\substack{i,j \\ i \geq a \text{ or } j \geq b}} A_{i,j} \geq N(1 -\alpha \beta + \epsilon)
\end{equation}
Now, it is basic analysis to show in a first step that the quadratic form 
\begin{equation}\label{eq:eq4}
\sum_{i,j=0}^{\sqrt{N}-1} \left( (A_{i,j} + B_{i,j})^2 -B_{i,j}^2 \right)
\end{equation}
attains its minimum under conditions (\ref{eq:eq1}), (\ref{eq:eq2}) and (\ref{eq:eq3}) if in (\ref{eq:eq2}) and (\ref{eq:eq3}) we have equality. In a second step, it can be shown that (\ref{eq:eq4}) attains its minimum under conditions (\ref{eq:eq1}), (\ref{eq:eq2}) and (\ref{eq:eq3}) (with equality sign) if all $A_{i,j}$ occurring in the sum of (\ref{eq:eq2}) and (\ref{eq:eq3}) have the same value. This means the minimum is attained if
\begin{equation*}
A_{i,j} =
\begin{cases}
& \frac{N(\alpha \beta - \epsilon)}{ab} = 1 -\frac{\epsilon}{\alpha \beta}, \ \text{ if } 0 \leq i < a \text{ and } 0 \leq j < b, \\
& \frac{N(1 -\alpha \beta + \epsilon)}{N - ab} = 1 + \frac{\epsilon}{1- \alpha \beta}, \ \text{ otherwise.}
\end{cases}
\end{equation*}
Note that each $B_{i,j}$ consists of $\frac{(2s-1)^2-1}{2} = 2s(s-1)$ summands $A_{x,y}$ and therefore, 
\begin{align*}
& NF_N^{(d)}(s) \geq \sum_{i=0}^{a-1} \sum_{j=0}^{b-1} \left( (A_{i,j} +B_{i,j})^2 -B_{i,j}^2 \right) +  \sum_{\substack{i,j \\ i \geq a \text{ or } j \geq b}} \left( (A_{i,j} + B_{i,j})^2 -B_{i,j}^2 \right) -N \\
& \geq ab \left(\left(\left(2s(s-1)+1\right)\left(1 - \frac{\epsilon}{\alpha \beta}\right)\right)^2 - \left(\left(2s(s-1)\right)\left(1-\frac{\epsilon}{\alpha \beta}\right)\right)^2\right) + \\
& + (N-ab) \left(\left(\left(2s(s-1)+1\right)\left(1 + \frac{\epsilon}{1-\alpha \beta}\right)\right)^2 - \left(\left(2s(s-1)\right)\left(1 + \frac{\epsilon}{1-\alpha \beta}\right)\right)^2\right) -N \\
& = (4s(s-1) +1)\left[ab \left(1- \frac{\epsilon}{\alpha \beta} \right)^2 +(N-ab) \left(1+ \frac{\epsilon}{1- \alpha \beta} \right)^2 \right] -N \\
& = N \left[ (4s(s-1)+1)\left(\lambda \left(1-\frac{\epsilon}{\lambda}\right)^2 +(1-\lambda)\left(1 + \frac{\epsilon}{1-\lambda}\right)^2 \right) -1 \right] \\
&=: NR_{\epsilon, \lambda}(s),  
\end{align*}
where $\lambda:= \alpha \beta$. By assumption, we have $ \lim_{N \to \infty} F_N^{(d)}(s) =(2s)^2$ for all positive integers $s$. Therefore, in order to derive a contradiction, it suffices to show that there exists an integer $s$ such that 
\begin{equation}\label{eq:eq5}
R_{\epsilon,\lambda}(s) > (2s)^2.
\end{equation}
The expression $R_{\epsilon, \lambda}(s) -(2s)^2$ can be viewed as a quadratic polynomial in $s$ with leading coefficient 
\begin{equation*}
4 \lambda \left(1-\frac{\epsilon}{\lambda}\right)^2 +(1-\lambda)\left(1 + \frac{\epsilon}{1-\lambda}\right)^2  -4 = 4 \frac{\epsilon^2}{\lambda(1 - \lambda)} >0.
\end{equation*}
Hence (\ref{eq:eq5}) holds for all $s$ large enough in dependence on $\epsilon$ and $\lambda$. \hfill $\square$
\section{Proof of Theorem \ref{TH:TH3}}
We again prove the result for $d=2$ only, as the general case is carried out quite analogously. 

There exists a constant $\rho$ with $0 < \rho < 1$ having the following property: For every pair $(\alpha_1, \alpha_2)$ there exist infinitely many $q \in \NN$ such that 
\begin{equation*}
\max(\lbrace q \alpha_1 \rbrace, \lbrace q \alpha_2 \rbrace) < \frac{\rho}{q^{1/2}},
\end{equation*}
see, e.g., \cite{not16}. Consider now such a $q$ and set 
\begin{equation*}
\frac{\theta}{q^{1/2}}=\max(\lbrace q \alpha_1 \rbrace, \lbrace q \alpha_2 \rbrace),
\end{equation*}
where $0 < \theta \leq \rho$. Let $A = A(q)$ be the minimal integer such that 
\begin{equation*}
\left( \left( \frac{1}{A \theta} \right)^{2/3} + 1 \right)^3 \theta^2 < \frac{1+\theta^2}{2}
\end{equation*}
holds, which is possible due to $\theta \leq \rho < 1$. Note that $A$ is the larger, the larger $\theta$ is. Hence, the values of $A$ are bounded by the value obtained for $\theta = \rho$. This $A$ will be denoted by $A_{\rho}$. Further, for the choice of $B=B(q):= \frac{2}{1+\theta^2}$, we have $B \geq \frac{2}{1+ \rho^2} >1$. Choose $L=L(q):= \left \lceil \left(\frac{1}{A \theta}\right)^{2/3} \right \rceil$ and the real $\tilde{\nu}$ such that
\begin{equation}\label{KR:eq1}
L \frac{\theta}{q^{1/2}} = \frac{1}{(A^2Lq-\tilde{\nu})^{1/2}},
\end{equation}
i.e., 
\begin{equation*}
\tilde{\nu} = A^2 Lq - \frac{q}{L^2 \theta^2}.
\end{equation*}
In the sequel, we will show that
\begin{equation}\label{KR:eq2}
\frac{q}{L^2\theta^2} \geq BLq,
\end{equation}
and consequently, 
\begin{equation}\label{KR:eq3}
\tilde{\nu} \leq (A^2-B)Lq.
\end{equation}
Clearly, equation (\ref{KR:eq2}) is equivalent to $BL^3\theta^2 \leq 1$. Now 
\begin{equation*}
BL^3 \theta^2 \leq B \left( \left( \frac{1}{A \theta} \right)^{2/3} +1 \right)^3 \theta^2 
\end{equation*}
and hence, (\ref{KR:eq2}) holds if 
\begin{equation*}
\left( \left( \frac{1}{A \theta}\right)^{2/3} +1 \right)^3 \theta^2 < \frac{1}{B} = \frac{1 + \theta^2}{2},
\end{equation*}
which is true due to the definition of $A$ and $B$, respectively. Let $\nu := \lfloor \tilde{\nu} \rfloor$ and $N:=A^2Lq - \nu$. Note, that by (\ref{KR:eq3}) we have  $N \ge  BLq$  and this, by the definition of $B$ and $L$ tends to infinity for $q$ to infinity.

We consider now the sequence elements $\boldsymbol{x}_1, \ldots, \boldsymbol{x}_{N=A^2Lq-\nu}$, and study the distances of the pairs 
\begin{equation*}
(\boldsymbol{x}_1,\boldsymbol{x}_{1 +qL}), \ (\boldsymbol{x}_2,\boldsymbol{x}_{2+qL}), \ \ldots, \ (\boldsymbol{x}_{N-qL},\boldsymbol{x}_N).
\end{equation*} 
Using the estimate (\ref{KR:eq3}), we derive that there are 
\begin{align}\label{KR:eq4}
N-qL =((A^2-1)L)q-\nu &\geq (B-1)Lq \nonumber
 \\
&\geq \left( \frac{2}{1+\rho^2} -1 \right)Lq \nonumber
 \\
& \geq A^2Lq \left( \frac{2}{1+ \rho^2} - 1  \right) \frac{1}{A^2_{\rho}} \nonumber
 \\
& \geq \gamma N 
\end{align}
such pairs, where $\gamma$ is a positive fixed constant independent on $q$. 

Now, for each such pair, we get (for $k=1, \ldots, N-qL$)
\begin{align*}
\| \boldsymbol{x}_k - \boldsymbol{x}_{k +qL} \|_{\infty} &= \max(\lbrace Lq \alpha_1 \rbrace, \lbrace Lq \alpha_2 \rbrace) \\
&=L\frac{\theta}{q^{1/2}}.
\end{align*}
Note that the second equality is true due to the inequality
\begin{equation*}
\frac{1}{\sqrt{N}} \leq L\frac{\theta}{q^{1/2}} \leq \frac{3}{\sqrt{N}},
\end{equation*}
which we will prove subsequently. 

First, note that we have
\begin{align*}
0 &< \sqrt{N} \left(\frac{1}{(A^2Lq-\tilde{\nu})^{1/2}} - \frac{1}{(A^2Lq- \nu )^{1/2}} \right) \\
& = \left( \frac{A^2Lq-\nu}{A^2Lq-\tilde{\nu}} \right)^{1/2} -1 = \left(1 + \frac{\tilde{\nu}- \nu}{(A^2Lq-\tilde{\nu})} \right)^{1/2}-1  \\
& < \left(1 + \frac{1}{BLq} \right)^{1/2} \leq \left(1 + \frac{1}{q} \right)^{1/2} <2,
\end{align*}
where we used the estimate (\ref{KR:eq3}) for $\tilde{\nu}$. Further, due to above estimate and (\ref{KR:eq1}), we also have
\begin{align*}
0 &< \sqrt{N} \left(L \frac{\theta}{q^{1/2}} - \frac{1}{\sqrt{N}} \right) \\
&= \sqrt{N} \left(\frac{1}{(A^2Lq-\tilde{\nu})^{1/2}} -\frac{1}{(A^2Lq-\nu)^{1/2}} \right) \\
&<2, 
\end{align*}
i.e., 
\begin{equation*}
\frac{1}{\sqrt{N}} \leq L\frac{\theta}{q^{1/2}} \leq \frac{3}{\sqrt{N}}.
\end{equation*}
This is valid for infinitely many $q$ and as a consequence thereof, there exists an $a \in \left\lbrace 0,1,\ldots, \left\lceil 3 \frac{100}{\gamma} \right\rceil \right\rbrace$ such that
\begin{equation*}
\frac{1+a \frac{\gamma}{100}}{\sqrt{N}} \leq L \frac{\theta}{q^{1/2}} \leq \frac{1+(a+1) \frac{\gamma}{100}}{\sqrt{N}} 
\end{equation*}
for infinitely many $q$. 

In the following, we consider a sequence $(q_l)$ of such $q$ with corresponding $(N_l)$ and define
\begin{align*}
s_1&:= 1 + a \frac{\gamma}{100} \\
s_2&:= 1 + (a+1) \frac{\gamma}{100}.
\end{align*}
Assume that $(\boldsymbol{x}_n)_{n \in \NN}$ were Poissonian, then we had
\begin{align*}
\lim_{l \to \infty} \Gamma_{i}(l):&= \\
& \lim_{l \to \infty} \frac{1}{N_l} \# \left\lbrace 1 \leq l \neq  m \leq N_i : \| \boldsymbol{x}_l - \boldsymbol{x}_m \|_{\infty} \leq \frac{s_i}{N^{1/2}} \right\rbrace = 4s_i^2,
\end{align*}
for $i=1,2$. 

But $ \Gamma_2(l) \geq \Gamma_1(l) + \gamma$ due to (\ref{KR:eq4}), which gives a contradiction as we also have 
\begin{equation*}
4s_2^2 - 4 s_1^2 = 4(s_2 +s_1)(s_2-s_1) \leq 32 \frac{\gamma}{100}.
\end{equation*}
\hfill $\square$
\section{Proof of Theorem \ref{TH:THMet}}
We adapt the steps of Lemma 3 of \cite{not5}, which we will shortly repeat here. \\ \\
In the sequel let $d \geq 2$ and $1,\alpha_1, \ldots, \alpha_d$ be linearly independent over the rationals. For $\boldsymbol{\alpha}=(\alpha_1, \ldots, \alpha_d)$, we denote by $\boldsymbol{I}_{N,s}(\boldsymbol{\alpha})$, for a fixed $s \asymp 1$ and $N$, with $(2s)^d \leq N$, the indicator function of the interval
\begin{equation*}
\left[-s/N^{1/d},s/N^{1/d} \right)^d
\end{equation*}
extended with period $\boldsymbol{1}$. Therefore, we can write the pair correlation function $F_N^{(d)}(s)$ (subsequently interpreted as a function of $\boldsymbol{\alpha}$) for the sequence $(\lbrace a_n \boldsymbol{\alpha} \rbrace)_{n \in \NN}$ as
\begin{equation*}
F_{N,s}^{(d)}(\boldsymbol{\alpha}):= \frac{1}{N} \sum_{\substack{1 \leq k, l \leq N \\ k \neq l}} \boldsymbol{I}_{N,s}(\boldsymbol{\alpha} (a_k -a_l)).
\end{equation*}
We consider the Fourier series expansion of $\boldsymbol{I}_{N,s}(\boldsymbol{\alpha})$, i.e., 
\begin{equation*}
\boldsymbol{I}_{N,s}(\boldsymbol{\alpha}) \sim \sum_{\boldsymbol{r} \in \ZZ^{d}} c_{\boldsymbol{r}} e(\left\langle \boldsymbol{r}, \boldsymbol{\alpha} \right\rangle), 
\end{equation*}
with
\begin{equation*}
c_{\boldsymbol{r}}= \int_{-s/N^{1/d}}^{s/N^{1/d}} \ldots \int_{-s/N^{1/d}}^{s/N^{1/d}} e^{-2 \pi i \sum_{i=1}^{d} r_i \alpha_i} \ d \alpha_1 \ldots d\alpha_d,
\end{equation*}
where $\boldsymbol{r} = (r_1, \ldots, r_d)$. Hence, we can write $c_{\boldsymbol{r}}=c_{r_1} \ldots c_{r_d}$, where $c_{r_j}=\int_{-s/N^{1/d}}^{s/N^{1/d}}e^{-2 \pi i r_j \alpha_j} \ d \alpha_j $, for $j=1, \ldots, d$ . Note that we have
\begin{equation}\label{eq:eqFour}
|c_{r_j}| \leq \min \left(\frac{2s}{N^{1/d}}, \frac{1}{|r_j|} \right).
\end{equation}
This gives the following estimate for the variance of $F_{N,s}^{(d)}(\boldsymbol{\alpha})$, where, due to the assumption on $s$, the constants implied by "$\ll$" are independent of $s$ (and of course independent of $N$).  
\begin{align*}
&\int_{[0,1)^d} \left( F_{N,s}^{(d)}(\boldsymbol{\alpha}) - \frac{(2s)^d(N-1)}{N} \right)^2 \  d \boldsymbol{\alpha} \\
&= \frac{1}{N^2} \int_{[0,1)^d} \left(\sum_{\substack{1 \leq k, l \leq N \\ k \neq l}}  \sum_{\boldsymbol{r} \in \ZZ^{d}\setminus{\lbrace\boldsymbol{0}\rbrace}} c_{\boldsymbol{r}} e(\left\langle \boldsymbol{r}, \boldsymbol{\alpha} (a_k -a_l)\right\rangle)  \right)^2 \  d \boldsymbol{\alpha} \\
&= \frac{1}{N^2} \sum_{\substack{1 \leq k, l, m, n  \leq N \\ k \neq l, m \neq n}}  \sum_{\boldsymbol{r}_1, \boldsymbol{r}_2  \in \ZZ^{d}\setminus{\lbrace\boldsymbol{0}\rbrace}} | c_{\boldsymbol{r}_1} c_{\boldsymbol{r}_2} | \times \\
& \times \int_{[0,1)^d} e(\left\langle \boldsymbol{r}_1, \boldsymbol{\alpha} (a_k -a_l)\right\rangle - \left\langle \boldsymbol{r}_2, \boldsymbol{\alpha} (a_m -a_n)\right\rangle ) \  d \boldsymbol{\alpha} \\
&= \frac{1}{N^2} \sum_{\substack{1 \leq k, l, m, n  \leq N \\ k \neq l, m \neq n}}  \sum_{\boldsymbol{r}_1, \boldsymbol{r}_2  \in \ZZ^{d}\setminus{\lbrace\boldsymbol{0}\rbrace}} | c_{\boldsymbol{r}_1} c_{\boldsymbol{r}_2} | \int_{[0,1)} e \left( \alpha_1 \left(r_1^{(1)} (a_k -a_l) - r_2^{(1)}(a_m-a_n) \right) \right) \  d \alpha_1 \ldots \times \\ 
& \times \int_{[0,1)} e \left( \alpha_d \left(r_1^{(d)} (a_k -a_l) - r_2^{(d)}(a_m-a_n) \right) \right) \  d \alpha_d,
\end{align*} 
with $\boldsymbol{r}_i=(r_i^{(1)}, \ldots, r_i^{(d)})$ for $i=1,2$. \\ \\
Introducing the representation function $r_N(v)$ of an integer $v$ defined as
\begin{equation*}
r_N(v):=  \# \lbrace 1 \leq k \neq l \leq N :  a_k -a_l =v \rbrace,
\end{equation*}
we can write the previous expression as 
\begin{equation}\label{eq:eqGCD}
\frac{1}{N^2} \sum_{v,w \in \ZZ \setminus{\lbrace 0 \rbrace}} r_N(v) r_N(w) \sum_{\substack{\boldsymbol{r}_1, \boldsymbol{r}_2  \in \ZZ^{d}\setminus{\lbrace\boldsymbol{0}\rbrace} \\ r_1^{(i)} v = r_2^{(i)}w }}' | c_{\boldsymbol{r}_1} c_{\boldsymbol{r}_2} |,
\end{equation}
where in the sum $\sum{'}$ the equality $r_1^{(i)} v = r_2^{(i)}w$ only needs to hold for these indices $i$ for which $r_1^{(i)}, r_2^{(i)} \neq 0$.
To furhter estimate the expression (\ref{eq:eqGCD}), we first derive the estimate (see \cite{not5} for details)
\begin{equation*}
\sum_{\substack{\boldsymbol{r}_1, \boldsymbol{r}_2  \in \ZZ^{d}\setminus{\lbrace\boldsymbol{0}\rbrace} \\ r_1^{(i)} v = r_2^{(i)}w }} | c_{\boldsymbol{r}_1} c_{\boldsymbol{r}_2}| \ll \frac{s^d \log N}{N} \frac{\gcd(v,w)}{\sqrt{|vw|}}, \qquad v,w \neq 0,
\end{equation*}
where the implied constant depends on $d$ and all entries of $\boldsymbol{r}_1$ and $\boldsymbol{r}_2$ are assumed to be non-zero. To see this, we recall that $r_1^{(i)} v = r_2^{(i)}w$ for $i=1, \ldots, d$ if and only if 
\begin{equation*}
r_1^{(i)} =\frac{h_iw}{\gcd(v,w)}, \text{ and } r_2^{(i)} =\frac{h_iv}{\gcd(v,w)} \quad \text{ for } i=1, \ldots, d,
\end{equation*}
where the $h_i$'s are some integers. Then, case distinctions according to the size of the $h_i$'s, i.e.,
\begin{align*}
& |h_i| \leq \left(\frac{N \gcd(v,w)}{s^d \max(|v|,|w|)} \right)^{1/d}=:\text{max}_{h_i}, \\
& \left(\frac{N \gcd(v,w)}{s^d \max(|v|,|w|)} \right)^{1/d} \leq |h_i| \leq  \left(\frac{N \gcd(v,w)}{s^d \min(|v|,|w|)} \right)^{1/d}=:\text{min}_{h_i}, \\
& |h_i| \geq  \left(\frac{N \gcd(v,w)}{s^d \min(|v|,|w|)} \right)^{1/d}
\end{align*}
and recalling the bounds on the Fourier coefficients (\ref{eq:eqFour}) gives the following estimate. For fixed $v$, $w$, and $\mathcal{D}:= \lbrace 1, \ldots, d \rbrace$, we have
\begin{align*}
&\sum_{\substack{\boldsymbol{r}_1, \boldsymbol{r}_2  \in \ZZ^{d}\setminus{\lbrace\boldsymbol{0}\rbrace} \\ r_1^{(i)} v = r_2^{(i)}w }} | c_{\boldsymbol{r}_1} c_{\boldsymbol{r}_2}| = \sum_{\substack{\boldsymbol{r}_1, \boldsymbol{r}_2  \in \ZZ^{d}\setminus{\lbrace\boldsymbol{0}\rbrace} \\ r_1^{(i)} v = r_2^{(i)}w }} | c_{r_1^{(1)}} \ldots c_{r_1^{(d)}} c_{r_2^{(1)}} \ldots c_{r_2^{(d)}}| \\
&\ll \sum_{\substack{\vartheta_1, \vartheta_2 \subseteq \mathcal{D} \\ \vartheta_1 \cap \vartheta_2 =\emptyset}} \sum_{\substack{h_i \leq \text{max}_{h_i} \\ i \in \vartheta_1}} \frac{s^{2|\vartheta_1|}}{N^{2|\vartheta_1|/d}} \sum_{\substack{\text{max}_{h_i} \leq h_i \leq \text{min}_{h_i} \\ i \in \vartheta_2}} \frac{s^{|\vartheta_2|}}{N^{|\vartheta_2|/d}} \frac{\gcd(v,w)^{|\vartheta_2|}}{\max(|v|,|w|)^{|\vartheta_2|}\prod_{i \in \vartheta_2} h_i} \times \\ 
&\times \sum_{\substack{h_i \geq \text{min}_{h_i} \\ i \in \mathcal{D}\setminus (\vartheta_1 \cup \vartheta_2)}} \frac{\gcd(v,w)^{2(d- |\vartheta_1 \cup \vartheta_2|)}}{|vw|^{d- |\vartheta_1 \cup \vartheta_2|}\prod_{i \in \mathcal{D}\setminus (\vartheta_1 \cup \vartheta_2)} h_i^2} \\ 
&\ll \sum_{\substack{\vartheta_1, \vartheta_2 \subseteq \mathcal{D} \\ \vartheta_1 \cap \vartheta_2 =\emptyset}} \left( \frac{N \gcd(v,w)}{s^d \max(|v|,|w|)} \right)^{|\vartheta_1|/d} \frac{s^{2|\vartheta|_1}}{N^{2|\vartheta_1|/d}} \times \\
& \times \log N \frac{s^{|\vartheta_2|}}{N^{|\vartheta_2|/d}} \frac{\gcd(v,w)^{|\vartheta_2|}}{\max(|v|,|w|)^{|\vartheta_2|}}  \frac{\gcd(v,w)^{2(d- |\vartheta_1 \cup \vartheta_2|)}}{|vw|^{d-|\vartheta_1 \cup \vartheta_2|}} \left( \frac{s^d}{N} \frac{\min(|v|,|w|)}{\gcd(v,w)} \right)^{(d-|\vartheta_1 \cup \vartheta_2|)/d} \\ 
& \ll \frac{s^d \log N}{N} \frac{\gcd(v,w)}{\sqrt{|vw|}}.
\end{align*}
Consequently, we have for the variance of $F_{N,s}^{(d)}(\boldsymbol{\alpha})$, using Theorem 4 of \cite{not20},  
\begin{align*}
&\int_{[0,1)^d} \left( F_{N,s}^{(d)}(\boldsymbol{\alpha}) - \frac{(2s)^d(N-1)}{N} \right)^2 \  d \boldsymbol{\alpha} \\
&= \frac{1}{N^2} \sum_{v,w \in \ZZ \setminus{\lbrace 0 \rbrace}} r_N(v) r_N(w) \sum_{\substack{\boldsymbol{r}_1, \boldsymbol{r}_2  \in \ZZ^{d}\setminus{\lbrace\boldsymbol{0}\rbrace} \\ r_1^{(i)} v = r_2^{(i)}w }}' | c_{\boldsymbol{r}_1} c_{\boldsymbol{r}_2} |  \\
& \ll\frac{s^d \log N}{N^3} \sum_{v, w \in \ZZ\setminus\lbrace 0 \rbrace} r_N(v) r_N(w) \frac{\gcd(v,w)}{\sqrt{|vw|}} \\
& \ll s^d (\log N)^{\mathcal{O}(1)} \exp \left(\mathcal{O}\left(\log \frac{N^3}{E(A_N)} \log \log N \right)^{1/2}\right) \frac{E(A_N)}{N^3}.
\end{align*} 
Then, following the lines of the proof of Theorem 6 of \cite{not20} allows to deduce the claim. \hfill $\square$

\section{Proof of Theorem \ref{TH:THMAX}}

As the proof uses exactly the same steps, except some minor technical changes, 
as the one in the one-dimensional case, see \cite{not8}, we will omit a detailed illustration. 
To prove the result in the one-dimensional case, we mention that most of the arguments are 
based on the additive structure of the integer sequence $(a_n)_{n \in \NN}$, i.e., 
the claim that there are "many" difference vectors $\boldsymbol{u}$ for which we have 
``many'' pairs $(k,l)$ such that $\boldsymbol{u}=\boldsymbol{r}_k -\boldsymbol{r}_l$ 
(cf., \textit{Property 1} in the proof of the main result of \cite{not8}) does not change 
in the multi-dimensional setting. For \textit{Property 2} in \cite{not8}, one has to 
consider now an interval of the form 
$\left[\beta, \beta + \frac{L}{\gamma N^{1/d}}\right)^d$, 
for some constants $\beta, \gamma, L>0$. 
The remaining arguments following \textit{Property~2}, however, do not change either and 
also note that in our definition of the supremum-norm, we work with the distance modulo 1 
instead of the absolute value.

\textbf{Author's Addresses:} \\ 
Gerhard Larcher, Lisa Kaltenb{\"o}ck and Wolfgang Stockinger, Institut f\"ur Finanzmathematik und Angewandte Zahlentheorie, Johannes Kepler Universit\"at Linz, Altenbergerstra{\ss}e 69, A-4040 Linz, Austria. \\ \\
Aicke Hinrichs and Mario Ullrich, Institut f\"ur Analysis, Johannes Kepler Universit\"at Linz, Altenbergerstra{\ss}e 69, A-4040 Linz, Austria.\\ \\ 
Email: aicke.hinrichs(at)jku.at, lisa.kaltenboeck(at)jku.at, gerhard.larcher(at)jku.at, wolfgang.stockinger(at)jku.at, mario.ullrich(at)jku.at.   
\end{document}